\pdfoutput=1

\documentclass[
12pt,
a4paper,
final,      
]{article}

\usepackage{artmacs}

\makeatletter
\def\namedlabel#1#2{\begingroup
    \def\@currentlabel{\protect{\textnormal{#2}}}%
    \label{#1}\endgroup
}

\newcommand{\SimpConst}[1]{S({#1})} 
\newcommand{\MultConst}[1]{M({#1})} 

\newcommand{\addpol}{A}
\newcommand{\subadd}{S}

\newcommand{\Pall}[1]{P_{#1}^{\text{all}}}
\newcommand{\Pone}[1]{P_{#1}}

\newcommand{\Ione}[1]{I_{#1}}

\newcommand{\Rone}[1]{R_{#1}}

\newcommand{\Qone}[1]{Q_{#1}}

\newcommand{\Sone}[1]{S_{#1}}

\newcommand{\Eone}[1]{E_{#1}}

\newcommand{\Aone}[1]{A_{#1}}

\newcommand{\qvar}{\mathbf{q}}
\newcommand{\degq}{\deg_{\qvar}}

\newcommand{\Qz}{\QQ \bbracket{z}}
\newcommand{\Qq}{\QQ (\qvar)}
\newcommand{\Qqz}{\Qq \bbracket{z}}

\newcommand{\combclass}[1]{{\mathcal #1}_{r}}
\newcommand{\ccP}{\combclass{P}}
\newcommand{\ccI}{\combclass{I}}
\newcommand{\ccR}{\combclass{R}}
\newcommand{\ccS}{\mathcal{S}_{r,s}}
\newcommand{\ccQ}{\mathcal{Q}_{r,s}}

\newcommand{\ngenfun}[1]{{\mathrm #1}_{r}} 
\newcommand{\nP}{\ngenfun{P}}
\newcommand{\nI}{\ngenfun{I}}
\newcommand{\nR}{\ngenfun{R}}
\newcommand{\nS}{\mathrm{S}_{r,s}}
\newcommand{\nQ}{\mathrm{Q}_{r,s}}

\newcommand{\nA}{\ngenfun{A}}

\newcommand{\sgenfun}[1]{{\mathsf #1}_{r}} 
\newcommand{\sP}{\sgenfun{P}}

\newcommand{\sI}{\sgenfun{I}}
\newcommand{\sR}{\sgenfun{R}}
\newcommand{\sS}{\mathsf{S}_{r,s}}
\newcommand{\sQ}{\mathsf{Q}_{r,s}}

\newcommand{\ncoeff}[1]{{\mathrm #1}} 
\newcommand{\ncP}[1]{\ncoeff{P}_{r,#1}}
\newcommand{\ncI}[1]{\ncoeff{I}_{r,#1}}
\newcommand{\ncR}[1]{\ncoeff{R}_{r,#1}}
\newcommand{\ncQ}[1]{\ncoeff{Q}_{r,#1,s}}
\newcommand{\ncE}[1]{\ncoeff{E}_{r,#1}}

\newcommand{\scoeff}[1]{{\mathsf #1}} 
\newcommand{\scP}[1]{\scoeff{P}_{r,#1}}
\newcommand{\scR}[1]{\scoeff{R}_{r,#1}}
\newcommand{\scE}[1]{\scoeff{E}_{r,#1}}
\newcommand{\scQ}[1]{\scoeff{Q}_{r,#1,s}}

\newcommand{\spaceR}[1]{R_{#1}^{*}}    
\newcommand{\spaceE}[1]{E_{#1}^{*}}

\begin{document}

\title{Survey on counting \\ special types of polynomials}
\pdftitle{Survey on counting special types of polynomials}
\pdfsubject{Survey on counting special multivariate and
  decomposable univariate polynomials over finite fields}
\author{Joachim von~zur~Gathen \& Konstantin Ziegler\\
B-IT, Universit\"at Bonn\\
D-53113 Bonn, Germany\\
\email{{gathen, zieglerk}@bit.uni-bonn.de}\\
\url{http://cosec.bit.uni-bonn.de/}
}
\pdfauthor{Joachim von zur Gathen, Konstantin Ziegler}

\maketitle

\begin{abstract}
Most integers are composite and most univariate polynomials over a
finite field are reducible. The Prime Number Theorem and a classical
result of Gau\ss\ count the remaining ones, approximately and exactly.

For polynomials in two or more variables, the situation changes dramatically.  Most
multivariate polynomials are irreducible.  This survey presents counting results
for some special classes of multivariate polynomials over a finite
field, namely the the reducible ones, the $s$-powerful ones (divisible
by the $s$-th power of a nonconstant polynomial), the relatively
irreducible ones (irreducible but reducible over an extension field),
the decomposable ones, and also for reducible space curves.
These come as exact formulas and as approximations with relative errors that
essentially decrease exponentially in the input size.

Furthermore, a univariate polynomial $f$ is decomposable if $f = g
\circ h$ for some nonlinear polynomials $g$ and $h$.  It is
intuitively clear that the decomposable polynomials form a small
minority among all polynomials.  The tame case, where the
characteristic $p$ of $\Fq$ does not divide $n = \deg f$, is fairly
well-understood, and we obtain closely matching upper and lower bounds
on the number of decomposable polynomials. In the wild case, where $p$
does divide $n$, the bounds are less satisfactory, in particular when $p$ is the smallest prime
divisor of $n$ and divides $n$ exactly twice.  The crux of the matter
is to count the number of collisions, where essentially different $(g, h)$ yield
the same $f$.  We present a classification of all collisions at degree
$n = p^{2}$ which yields an exact count of those decomposable
polynomials.
\end{abstract}

\begin{keywords}
  counting special polynomials, finite fields,  combinatorics on polynomials,
  generating functions,   analytic combinatorics,   asymptotic behavior,
 multivariate polynomials, polynomial decomposition, Ritt's Second Theorem
\end{keywords}
\pdfkeywords{counting special polynomials, finite fields,  combinatorics on polynomials,
  generating functions,   analytic combinatorics,   asymptotic behavior,
 multivariate polynomials, polynomial decomposition, Ritt's Second Theorem}

\begin{AMS}
00B25, 11T06, 12Y05
\end{AMS}

\section{Introduction}

Most integers are composite and most univariate
polynomials over a finite field are reducible.  The classical results
of the Prime Number Theorem and a theorem of Gau\ss\ present
approximations saying that
randomly chosen integers up to $x$ or polynomials of degree
up to $n$ are prime or irreducible with probability about $1/\ln x$
or $1/n$, respectively.

Concerning special classes of univariate polynomials over a finite
field, \cite{zsi94} counts those with a given number of distinct roots
or without irreducible factors of a given degree.  In the same
situation, \cite{art24} counts the irreducible ones in an arithmetic
progression and \cite{hay65} generalizes
these results.  \cite{coh69a} and \cite{car87}
count polynomials with certain factorization patterns and \cite{wil69}
those with irreducible factors of given degree.  Polynomials that
occur as a norm in field extensions are studied by \cite{goglut81}.

In two or more variables, the situation changes dramatically. Most
multivariate polynomials are irreducible.  \cite{car63a} provides the
first count of irreducible multivariate polynomials.  In \cite{car65},
he goes on to study the fraction of irreducibles when bounds on the
degrees in each variable are prescribed; see also \cite{coh68}.  In
this survey, we opt for bounding the total degree because it has the
charm of being invariant under invertible linear transformations.
\cite{gaolau02} consider the counting problem in yet another model, namely
where one variable occurs with maximal degree.  The natural generating
function (or zeta function) for the irreducible polynomials in two or
more variables does not converge anywhere outside of the
origin. \cite{wan92} notes that this explains the lack of a simple
combinatorial formula for the number of irreducible polynomials. But
he gives a $p$-adic formula, and also a (somewhat complicated)
combinatorial formula.  For further references, see
\citet[Section~3.6]{mulpan13}.

In the bivariate case, \cite{gat08-incl-gat07}
proves precise approximations with an exponentially decreasing relative
error.  \Citet*{gatvio13} extend those results to multivariate
polynomials and give further information such as exact formulas and
generating functions.   \cite{bod08} gives a recursive formula for the number of irreducible
bivariate polynomials and remarks on a generalization for more than
two variables; he follows up with \cite{bod10}.

We present exact formulas for the numbers of reducible
(Sections~\ref{sec:gen}-\ref{sec:red}), $s$-powerful (Section~\ref{sec:powerful}),
and relatively irreducible polynomials (Section~\ref{sec:rel_irr}).  The
formulas also yield simple, yet precise, approximations to these numbers, with rapidly decaying
relative errors.

Geometrically, a single polynomial corresponds to a hypersurface, that
is, to a cycle in affine or projective space of codimension 1.  This
correspondence preserves the respective notions of reducibility.  Thus,
Sections~\ref{sec:gen}-\ref{sec:red} can also be viewed as counting reducible
hypersurfaces, in particular, planar curves, and
Section~\ref{sec:powerful} those with an $s$-fold component.  From a
geometric perspective, these results say that almost all hypersurfaces
are irreducible.  Can we say something
similar for other types of varieties?
\cite{cesgat13} give an affirmative answer for curves in $\PP^{r}$ for arbitrary $r$. A first question is how
to parametrize the curves. Moduli spaces only include irreducible curves, and systems of defining
equations do not work except for complete intersections. The natural parametrization is by the Chow
variety $C_{r,n}$ of curves of degree $n$ in $\PP^{r}$, for some
fixed $r$ and $n$. The foundation of this approach is
a result by \cite{eishar92}, who identified the irreducible components of $C_{r,n}$ of maximal
dimension.  We present the counting results in Section~\ref{sec:reducible-curves}.

It is intuitively clear that the decomposable polynomials form a small
minority among all multivariate polynomials over a
field.  \Citet{gat10a} gives a quantitative version of this
intuition (see Section~\ref{sec:multi-decomp}).
The number of multivariate decomposable polynomials is also studied by
\cite{boddeb09}.

This concludes the first half (\autoref{sec:exact}) of our survey, dealing with multivariate
polynomials.   The second half (\autoref{sec:univ-decomp}) is devoted to
counting univariate decomposable polynomials.

Some of the results in this survey are from joint work with Raoul
Blankertz, Eda Cesaratto, Mark Giesbrecht, Guillermo Matera, and
Alfredo Viola.

A version of this paper is to appear in \cite{RICAM2013}. The final
publication will be available at Springer after publication.

\section{Counting multivariate polynomials}
\label{sec:exact}

We work in the polynomial ring $F[x_{1}, \dots ,x_{r}]$  in
$r\geq 1$ variables over a field $F$ and consider polynomials
with total degree equal to some nonnegative integer~$n$:
\[
\Pall{r,n}(F) = \{ f \in F[x_{1},\dots,x_{r}] \colon \deg f = n\}.
\]
The polynomials of degree at most $n$
form an $F$-vector space of dimension~$\binom{r+n}{r}$.

The property of a certain polynomial to be reducible, squareful,
relatively irreducible, or decomposable is shared with all polynomials associated to
the given one. For counting them, it is sufficient to
take one representative. We choose an arbitrary monomial order, say, the degree-lexicographic one,
so that the monic polynomials are those with leading coefficient 1, and write
\[
\Pone{r,n}(F) = \{f \in \Pall{r,n}(F) \colon f \text{ is monic}\}.
\]

We use two different methodologies to obtain such bounds: generating
functions and combinatorial counting.  The usual approach, see
\cite{flased09}, of analytic combinatorics on series with integer
coefficients leads, in our case, to power series that diverge
everywhere (except at $0$).  We have not found a way to make this
work.  Instead, we use power series with symbolic coefficients, namely
rational functions in a variable representing the field size.
Several useful relations from standard analytic combinatorics carry
over to this new scenario.  In a first step, this yields in a
straightforward manner an exact formula for the number under
consideration (\autoref{pro:R_exact_by_recursion}).  This formula is, however, not very
transparent.  Even the leading term is not immediately visible.

In a second step, coefficient comparisons yield easy-to-use
approximations to our number (\autoref{thm:R-from-gen}).  The relative error is
exponentially decreasing in the bit size of the data.  Thus,
\autoref{thm:R-from-gen} gives a ``third order'' approximation for the number of reducible
polynomials, and thus a ``fourth order'' approximation for the
irreducible ones.  The error term is in the big-Oh form and thus
contains an unspecified constant.

In a third step, a different method, namely some combinatorial
counting, yields ``second order'' approximations with explicit
constants in the error term (\autoref{thm:red}).

The results of Sections~\ref{sec:gen}-\ref{sec:rel_irr} are from \citet*{gatvio13} unless otherwise
attributed, those of Section~\ref{sec:reducible-curves} are from
\cite*{cesgat13}, and those of Section~\ref{sec:multi-decomp} are from \citet{gat10a}.

\subsection{Exact formula for reducible polynomials}
\label{sec:gen}

To study reducible polynomials, we consider the following subsets of $\Pone{r,n}(F)$:
\begin{align}
  \Ione{r,n}(F) & = \{ f \in \Pone{r,n}(F) \colon f \text{ is irreducible} \}, \\
  \Rone{r,n}(F) & = \Pone{r,n}(F) \mysetminus \Ione{r,n}(F).
\end{align}
In the usual notions, the polynomial $1$ is neither reducible nor
irreducible.  In our context, it is natural to have $\Rone{r,0}(F) =
\{ 1 \}$ and $ \Ione{r,0}(F) = \varnothing$.

The sets of polynomials
\begin{align}
  \ccP & = \bigcup_{n \geq 0} \Pone{r,n} (\FF_{q}), \\
  \ccI & = \bigcup_{n \geq 0} \Ione{r,n} (\FF_{q}), \\
  \ccR & = \ccP \mysetminus \ccI,
\end{align}
are combinatorial classes with the total degree as size functions and
we denote the corresponding generating functions by $\nP, \nI, \nR \in
\ZZ_{\geq 0} \bbracket{z}$, respectively.  Their coefficients are
\begin{align}
\ncP{n} & =  \# \Pone{r,n} (\FF_{q}) = q^{\binom{r+n}{r}-1} \frac{1-q^{-\binom{r+n-1}{r-1}}}{1-q^{-1}}, \label{eq:44} \\
\ncR{n} & =  \#  \Rone{r,n} (\FF_{q}), \\
\ncI{n} & =  \#  \Ione{r,n} (\FF_{q}), \label{eq:13}
\end{align}
respectively, dropping the finite field $\FF_{q}$ with $q$ elements
from the notation.
By definition, $\ccP$ equals the disjoint union of $\ccR$ and $\ccI$, and therefore
\begin{equation}
\label{eq:42}
\nR = \nP - \nI.
\end{equation}
By unique factorization, every element in $\ccP$ corresponds to an
unordered finite sequence of elements in $\ccI$, where repetition is
allowed, and therefore
\begin{equation}
\nI  = \sum_{k \geq 1} \frac{\mu (k)}{k} \log \nP (z^{k}) \label{eq:28}
\end{equation}
by \citet[Theorem~I.5]{flased09}, where $\mu$ is the number-theoretic M\"obius-function.
A resulting algorithm is easy to program and returns exact results
with lightning speed.

This approach quickly leads to explicit formulas. A \emph{composition}
of a positive integer $n$ is a sequence $j = (j_{1}, j_{2}, \dots, j_{\abs{j}})$ of
positive integers $j_{1}, j_{2}, \dots, j_{\abs{j}}$ with $j_{1} +
j_{2} + \dots + j_{\abs{j}} = n$, where $\abs{j}$ denotes the length
of the sequence.  We define the set
\begin{equation}
  \label{eq:10}
M_{n}  = \{ \text{compositions of $n$} \}.
\end{equation}
  This standard
combinatorial notion is not to be confused with the composition of
polynomials, which we discuss in Sections~\ref{sec:multi-decomp} and \ref{sec:univ-decomp}.

\begin{theorem}[Exact counting]
  \label{pro:R_exact_by_recursion} Let $r \geq 1$, $q \geq 2$,
  $\ncP{n}$ as in \eqref{eq:44}, and $\ncI{n}$ the number of
  irreducible monic $r$-variate polynomials of degree $n$ over
  $\Fq$. Then we have
\begin{equation}
  \label{eq:12}
  \begin{split}
    \ncI{0} & = 0, \\
    \ncI{n} & = - \sum_{k \,\mid\, n} \frac{\mu (k)}{k} \sum_{j \in
      M_{n/k}} \frac{(-1)^{\abs{j}}}{\abs{j}} \ncP{j_{1}} \ncP{j_{2}} \cdots \ncP{j_{\abs{j}}},
  \end{split}
\end{equation}
for $n \geq 1$, and therefore for the number $\ncR{n}$ of reducible monic $r$-variate polynomials of degree $n$ over $\Fq$
  \begin{align}
    \ncR{0} & = 1, \\
    \ncR{n} & = \ncP{n} + \sum_{k \,\mid\, n} \frac{\mu (k)}{k} \sum_{j \in M_{n/k}} \frac{(-1)^{\abs{j}}}{\abs{j}} \ncP{j_{1}} \ncP{j_{2}} \cdots \ncP{j_{\abs{j}}},
  \end{align}
for $n \geq 1$.
\end{theorem}
The formula of \autoref{pro:R_exact_by_recursion} is exact but
somewhat cumbersome.  The following two sections provide simple yet
precise approximations, with rapidly decaying error terms.

\subsection{Symbolic approximation for reducible polynomials}
\label{sec:symb-appr-reduc}

For $r \geq 2$, the power series $\nP$, $\nI$, and $\nR$ do not converge anywhere except at 0, and the standard asymptotic arguments of analytic combinatorics are inapplicable.  We now deviate from this approach and move from power series in $\Qz$ to power series in $\Qqz$, where $\qvar$ is a symbolic variable representing the field size.  For $r \geq 2$ and $n \geq 0$ we let
\begin{equation}
  \label{eq:19}
\scP{n}(\qvar) = \qvar^{\binom{r+n}{r}-1}\frac{1-\qvar^{-\binom{r+n-1}{r-1}}}{1-\qvar^{-1}} \in \ZZ [\qvar]
\end{equation}
in analogy to \eqref{eq:44}. We define the power series $\sP, \sI, \sR \in \Qqz$ by
\begin{align}
  \sP (\qvar, z) & = \sum_{n \geq 0} \scP{n}(\qvar) z^{n}, \label{eq:5} \\
  \sI (\qvar, z) & = \sum_{k \geq 1} \frac{\mu (k)}{k} \log \sP (\qvar, z^{k}), \label{eq:relationIP} \\
  \sR (\qvar, z) & = \sP (\qvar, z) - \sI (\qvar, z). \label{eq:70}
\end{align}
Then $\scR{n}(\qvar)$ denotes the coefficient of $z^{n}$ in $\sR$ and
counts symbolically the reducible monic $r$-variate polynomials of degree $n$.

For nonzero $f \in \QQ(\qvar)$, $\degq f$ is the degree of $f$, that is, the numerator
degree minus the denominator degree.
The appearance of $O(\qvar^{-m})$ with a positive integer $m$ in an
equation means the existence of some $f$ with degree at most $-m$ that makes the
equation valid.  If a term $O(\qvar^{-m})$ appears, then we may
conclude a numerical
asymptotic result for growing prime powers $q$.

\begin{theorem}[Symbolic approximation]
\label{thm:R-from-gen}
  Let $r\geq 2$ and
  \begin{equation}
    \rho_{r,n}(\qvar)
= \qvar^{\binom{r+n-1}{r}+r-1}
    \frac{1-\qvar^{-r}}{(1-\qvar^{-1})^2}  \in \QQ(\qvar). \label{eq:rho}
\end{equation}
Then the symbolic formula $\scR{n}(\qvar)$ for the number of
reducible monic $r$-variate polynomials of degree $n$ over $\Fq$ satisfies
\begin{align}
\scR{0}(\qvar) & = 1, \quad \scR{1} (\qvar)= 0, \quad \scR{2} (\qvar)= \frac{ \rho_{r,2}(\qvar) }{2} \cdot (1-\qvar^{-r-1}), \\
\scR{3} (\qvar)& = \rho_{r,3}(\qvar) \Bigl( 1-\qvar^{-r(r+1)/2} +
\qvar^{-r(r-1)/2}
\frac{1-2\qvar^{-r}+2\qvar^{-2r-1}-\qvar^{-2r-2}}{3(1-\qvar^{-1})}\Bigr),
\\
\scR{4} (\qvar)& = \rho_{r,4} (\qvar) \cdot \Bigl( 1 +
\qvar^{-\binom{r+1}{3}} \cdot  \frac{1 +
  O(\qvar^{-r(r-1)/2})}{2 (1-\qvar^{-r})} \Bigr), \label{eq:48} \\
\intertext{and for $n \geq 5$}
\scR{n} (\qvar)& = \rho_{r,n} (\qvar) \cdot \Bigl( 1 + \qvar^{-\binom{r+n-2}{r-1}+r(r+1)/2} \cdot  \frac{1+ O(\qvar^{-r(r-1)/2})}{1-\qvar^{-r}} \Bigr). \label{eq:26}
\end{align}
\end{theorem}

\cite{ale06} lists $(\#
I_{r,n}(\Fq))_{n \geq 0}$ as A115457--A115472 in The On-Line Encyclopedia
of Integer Sequences, for $2 \leq r \leq 6$ and prime $q \leq 7$.  \citet[Theorem~7]{bod08} states (in our notation)
\[
1-\frac{\# I_{r, n}}{ \# P_{r, n} } \sim q^{-\binom{n+r-1}{r-1}-r}
\frac{1-q^{-r}}{1-q^{-1}}.
\]
\cite{houmul09} provide results for
$\# I_{r,n} (\FF_q)$.  These do not yield error bounds for the approximation
of $\# R_{r,n} (\FF_q)$.  \cite{bod10} also uses \eqref{eq:28} to
claim a result similar to \eqref{eq:26}.

\subsection{Explicit bounds for reducible polynomials}
\label{sec:red}

The third approach by ``combinatorial counting'' is somewhat more involved.  The payoff of this
additional effort is an explicit relative error bound.
However, the calculations are sufficiently complicated for us to stop at the
first error term.  Thus we replace the asymptotic $1 +
  O(\qvar^{-r(r-1)/2})$ in \eqref{eq:26} by $1/(1-q^{-1})$.

\begin{theorem}[Explicit approximation]
  \label{thm:red}
Let $r, q \geq 2$, and $\rho_{r,n}$ as in
  \autoref{thm:R-from-gen}. For the number $\#
  \Rone{r,n}(\Fq)$ of reducible monic $r$-variate polynomials of
  degree $n$ over $\Fq$ we have
  \begin{gather}
  \# \Rone{r,0} (\FF_{q}) = 1, \quad \# \Rone{r,1} (\FF_{q}) = 0,
  \quad \#
  \Rone{r,2} (\FF_q) = \frac{ \rho_{r,2}(q)}{2} \cdot (1-q^{-r-1}), \\
\begin{split}
    \left| \# \Rone{r,3} (\FF_q)- \rho_{r,3}(q) \right| &  = \rho_{r,3}(q) \cdot q^{-r(r-1)/2}
      \frac{1-2q^{-r}+2q^{-2r-1}-q^{-2r-2}}{3(1-q^{-1})} \\
& \leq
    \rho_{r,3}(q) \cdot q^{-r(r-1)/2}, \\
\intertext{and for $n \geq 4$}
  \left| \# \Rone{r,n} (\FF_q)-
      \rho_{r,n}(q) \right| & \leq \rho_{r,n}(q) \cdot
    \frac{q^{-\binom{r+n-2}{r-1}+r(r+1)/2}}{(1-q^{-1})(1-q^{-r})} \label{eq:16} \\
& \leq \rho_{r,n}(q) \cdot 3 q^{-\binom{r+n-2}{r-1}+r(r+1)/2}.
\end{split}
\end{gather}
\end{theorem}

\begin{remark}
  \label{rem:exp_decay}
  How close is our relative error estimate to being exponentially
  decaying in the input size? The usual dense representation of a
  polynomial in $r$ variables and of degree $n$ requires $b_{r,n} =
  \binom{r+n}{r}$ monomials, each of them equipped with a coefficient
  from $\FF_q$, using about $\log_2 q$ bits. Thus the total input size
  is about $\log_2 q \cdot b_{r,n}$ bits.  This differs from $\log_2 q
  \cdot (b_{r-1,n-1}-b_{r-1,2})$ by a factor of
  \[
  \frac{b_{r,n}}{b_{r-1,n-1}-b_{r-1,2}} <
  \frac{b_{r,n}}{\frac{1}{2}b_{r-1,n-1}} = \frac{2(n+r)(n+r-1)}{nr}.
  \]
  Up to this polynomial difference (in the exponent), the relative
  error is exponentially decaying in the bit size of the input, that
  is, $(\log q)$ times the number of coefficients in the usual dense
  representation.  In particular, it
  is exponentially decaying in any of the parameters $r$, $n$, and $\log_2
  q$, when the other two are fixed.
\end{remark}

\subsection{Powerful polynomials}
\label{sec:powerful}

For an integer $s \geq 2$, a polynomial is called
\emph{$s$-powerful} if it is divisible by the $s$th power of some
nonconstant polynomial, and \emph{$s$-powerfree} otherwise; it is
\emph{squarefree} if $s=2$. Let
\begin{align}
  \Qone{r,n,s}   (F) & = \{f\in \Pone{r,n}  (F) \colon f \text{ is
    $s$-powerful}\}, \\
  \Sone{r,n,s} (F) & = \Pone{r,n}(F) \mysetminus \Qone{r,n,s}(F).
\end{align}
As in the previous section, we restrict our attention to a finite field
$F= \FF_q$, which we omit from the notation.

For the approach by generating functions, we consider the
combinatorial classes $\ccQ = \bigcup_{n \geq 0} \Qone{r,n,s}$ and
$\ccS = \ccP \mysetminus \ccQ$. Any monic polynomial $f$ factors
uniquely as $f=g\cdot h^s$ where $g$ is a monic $s$-powerfree
polynomial and $h$ an arbitrary monic polynomial, hence
\begin{equation}
  \label{eq:61}
\nP = \nS \cdot \nP ( z^{s})
\end{equation}
and  by definition $\nQ = \nP
- \nS$ for the generating functions of $\ccS$ and $\ccQ$,
respectively.  For univariate polynomials, \cite{car32} derives \eqref{eq:61}
directly from generating functions to prove the counting formula
\eqref{eq:71} for $r=1$.  \citet[Section~1.1]{flagou01} use
\eqref{eq:61} for $s=2$ to count univariate squarefree polynomials,
see also
\citet[Note~I.66]{flased09}.

As in \autoref{pro:R_exact_by_recursion}, this approach quickly leads to explicit formulas.
\begin{theorem}[Exact counting]
\label{pro:Q_exact_by_recursion}
For $r \geq 1$, $q,s \geq 2$, $\ncP{n}$ as in \eqref{eq:44}, and
$M_{n}$ as in \eqref{eq:10}, we have for the number $\ncQ{n} = \#
\Qone{r,n,s} (\FF_{q})$ of $s$-powerful monic $r$-variate polynomials
of degree $n$ over $\Fq$
\begin{equation}
  \ncQ{n} = - \sum_{\substack{1 \leq i
      \leq n/s \\ j \in M_{i}}} (-1)^{\abs{j}}  \ncP{j_{1}}
  \ncP{j_{2}} \cdots \ncP{j_{\abs{j}}} \ncP{n-is} \label{eq:71}.
\end{equation}
\end{theorem}

To study the asymptotic behavior of $\ncQ{n}$ for $r \geq 2$ we again deviate from the standard approach and move to power series in $\Qqz$.  With $\sP$ from \eqref{eq:5}, we define $\sS, \sQ \in \Qqz$ by
  \begin{align}
    \sP & = \sS \cdot \sP(z^{s}), \label{def:S} \\
    \sQ & = \sP - \sS.
\end{align}

The approach by generating functions now yields the following result.  Its ``general'' case is \ref{item:13}.  We give exact expressions in special cases, namely for $ n < 3s$ in \ref{item:17} and for $(n,s)=(6,2)$ in \ref{item:18}, which also apply when we substitute the size $q$ of a finite field $\FF_{q}$ for $\qvar$.

\begin{theorem}[Symbolic approximation]
\label{thm:Q_by_gen}
  Let $r, s \geq 2$, $n\geq 0$, and
  \begin{align}
  \eta_{r,n,s} (\qvar) & = \qvar^{\binom{r+n-s}{r} + r-1}
  \frac{(1-\qvar^{-r}) (1-\qvar^{-\binom{r+n-s-1}{r-1}})}{(1-\qvar^{-1})^2} \in
  \QQ (\qvar), \label{eq:89} \\
  \delta & = \binom{r+n-s}{r}-\binom{r+n-2s}{r} - \frac{r(r+1)}{2}.
  \end{align}
  Then the symbolic formula $\scQ{n} (\qvar)$ for the number of $s$-powerful monic $r$-variate polynomials
of degree $n$ over $\Fq$ satisfies the following.
  \begin{ronumerate}
  \item \label{item:4} If $n \geq 2s$, then $\delta \geq r$.
  \item \label{item:17} \begin{equation}
      \label{eq:50}
      \scQ{n}(\qvar) = \begin{cases}
0 & \text{for $n<s$,}\\
\eta_{r,n,s}(\qvar) & \text{for $s \leq n<2s$,} \\
\eta_{r,n,s} (\qvar) \biggl( 1 + \qvar^{-\delta} \cdot \frac{1-\qvar^{-\binom{n+r-2s-1}{r-1}}}{1-\qvar^{-\binom{n+r-s-1}{r-1}}} & \\
\quad \cdot \Bigl( \frac{1-\qvar^{-r(r+1)/2}}{1-\qvar^{-r}} - \qvar^{-r(r-1)/2}\frac{1-\qvar^{-r}}{1-\qvar^{-1}}\Bigr)\biggr) & \text{for $2s \leq n<3s$.}
\end{cases}
\end{equation}

  \item \label{item:18} For $(n,s)=(6,2)$, we have
    \begin{align}
\scoeff{Q}_{r,6,2} (\qvar) 
& = \eta_{r,6,2} (\qvar) \bigl( 1 + \qvar^{-\delta+(r-2)(r-1)(r+3)/6}(1+O(\qvar^{-1}))\bigr). \label{eq:86}
    \end{align}
  \item\label{item:13} For $n \geq 2s$ and $(n,s)\neq (6,2)$, we have
    \begin{equation}
\label{eq:49}
      \scQ{n} (\qvar) = \eta_{r,n,s}(\qvar) \big(1 + \qvar^{-\delta}
      (1 + O(\qvar^{-1}))\big).
    \end{equation}
  \end{ronumerate}
\end{theorem}
For $r\geq 3$, we can replace $1+O(\qvar^{-1})$ in \eqref{eq:86} by
$\qvar^{-1}+O(\qvar^{-2})$.  The combinatorial approach replaces the asymptotic
$1+O(\qvar^{-1})$ for $n \geq 3s$ with an explicit bound.
For $n<3s$ the exact formula \eqref{eq:50} of \autoref{thm:Q_by_gen}~\ref{item:17} applies.

\begin{theorem}[Explicit approximation]
\label{thm:Q_by_map}
Let $r,s,q\geq 2$, $\# \Qone{r,n,s}(\Fq)$ the number
of $s$-powerful monic $r$-variate polynomials of degree $n$ over
$\Fq$, and $\eta_{r,n,s}$ and $\delta$ as in
\autoref{thm:Q_by_gen}.
\begin{ronumerate}
\item \label{item:8} For $(n,s) = (6,2)$, we have $\delta = r(r+1)(r^{2}+9r+2)/24$ and
\begin{equation}
\label{eq:53}
  \left| \# \Qone{r,6,2} (\FF_q) - \eta_{r,6,2}(q) \right| \leq
  \eta_{r,6,2}(q) \cdot 2 q^{-\delta+(r-2)(r-1)(r+3)/6}.
\end{equation}
  \item \label{item:9} For $n \geq 3s$ and $(n,s) \neq (6,2)$, we have
\begin{equation}
\label{eq:18}
  \left| \# \Qone{r,n,s} (\FF_q) - \eta_{r,n,s}(q) \right| \leq
  \eta_{r,n,s}(q) \cdot 6 q^{-\delta}.
\end{equation}
\end{ronumerate}
\end{theorem}

  As noted in \autoref{rem:exp_decay} for reducible polynomials, the
  relative error term is (essentially) exponentially decreasing in the
  input size, and exponentially decaying in any of the parameters $r$,
  $n$, $s$, and $\log_2 q$, when the other three are fixed.

\subsection{Relatively irreducible polynomials}
\label{sec:rel_irr}

A polynomial over $F$ is \emph{absolutely
  irreducible} if it is irreducible over an algebraic closure of $F$,
and \emph{relatively irreducible} (or \emph{exceptional}) if it
is irreducible over $F$ but factors over some extension field of $F$.  We
define
\begin{align}
  \Aone{r,n} (F) & = \{ f \in \Pone{r,n}(F) \colon f \text{ is absolutely irreducible}\} \subseteq \Ione{r,n}(F),\\
\Eone{r,n} (F) & = \Ione{r,n} (F) \mysetminus \Aone{r,n}(F) \label{eq:66}.
\end{align}
As before, we restrict ourselves to finite fields and recall that all our polynomials are monic.
We relate the generating function
$\nA(\FF_{q})$ of $\# \Aone{r,n} (\FF_{q})$ to the generating function $\nI (\FF_{q})$ of
irreducible polynomials as introduced in Section~\ref{sec:gen} and obtain
\begin{align}
  [z^{n}] \, \nI (\FF_{q}) & = \sum_{k \,\mid\, n } \frac{1}{k} \sum_{s
    \,\mid\, k} \mu (k/s) \cdot [z^{n/k}] \, \nA (\FF_{q^{s}}), \\
  [z^{n}] \, \nA (\FF_{q}) & = \sum_{k \,\mid\, n} \frac{1}{k} \sum_{s \,\mid\, k} \mu (s)
\cdot  [z^{n/k}] \, \nI (\FF_{q^{s}}) \label{eq:43}
\end{align}
with  M\"{o}bius inversion.  For an explicit formula, we combine the expression for $\ncI{n} (\FF_{q})$ from \autoref{pro:R_exact_by_recursion} with \eqref{eq:43}.

\begin{theorem}[Exact counting]
\label{pro:E_exact_by_recursion}
For $r,n \geq 1$, $q \geq 2$, $M_{n}$ as in \eqref{eq:10}, $\ncP{n}$
as in \eqref{eq:44}, and $\ncI{n}$ as in \eqref{eq:13}, we have for the
number $\ncE{n}$ of relatively irreducible monic $r$-variate
polynomials of degree $n$ over $\Fq$
  \begin{align}
\ncE{0} (\FF_{q})& = 0, \\
\ncE{n} (\FF_{q}) & =  -\sum_{1 < k \,\mid\, n} \frac{1}{k} \sum_{s \,\mid\, k} \mu (s) \ncI{n/k} (\FF_{q^{s}}) \label{eq:60} \\
& = \sum_{1< k \,\mid\, n} \frac{1}{k} \sum_{\substack{s \,\mid\, k \\
    m \,\mid\, n/k }} \frac{\mu (s) \mu (m)}{m} \\
& \quad \quad \quad \quad \cdot\sum_{j \in M_{n/(km)}} \frac{(-1)^{\abs{j}}}{\abs{j}} \ncP{j_{1}} (\FF_{q^{s}}) \ncP{j_{2}} (\FF_{q^{s}}) \cdots \ncP{j_{\abs{j}}} (\FF_{q^{s}}).
  \end{align}
\end{theorem}

The approach by generating functions gives the following result.

\begin{theorem}[Symbolic approximation]
\label{thm:E_by_gen}
Let $r, n \geq 2$, let $\ell$ be the smallest prime divisor of $n$, and
 \begin{align}
    \epsilon_{r,n}(\qvar) & = \frac{\qvar^{\ell ( \binom{r+n/\ell}{r} - 1)}}{\ell(1-\qvar^{-\ell})} \in \QQ(\qvar), \label{eq:epsilon} \\
\kappa & = (\ell-1) (\binom{r-1+n/\ell}{r-1}-r ) + 1.
 \end{align}
Then the symbolic formula $\scE{n}(\qvar)$ for the number of relatively irreducible monic $r$-variate
polynomials of degree $n$ over $\Fq$ satisfies the following.
\begin{ronumerate}
  \item \label{enum3:i} $\scE{1} (\qvar)= 0$.
  \item \label{enum3:ii} If $n$ is prime, then
    \begin{align}
    \scE{n} (\qvar) & = \epsilon_{r,n}(\qvar) (1-\qvar^{-nr}) \Bigl( 1 -
      \qvar^{-r(n-1)}\frac{(1-\qvar^{-r})(1-\qvar^{-n})}{(1-\qvar^{-1})(1-\qvar^{-nr})}\Bigr).
    \end{align}
    \item \label{enum3:iii} If $n$ is composite, then $\kappa \geq 2$ and
\[
    \scE{n} (\qvar) = \epsilon_{r,n}(\qvar) (
 1 + O (\qvar^{- \kappa}) ).
\]
\end{ronumerate}
\end{theorem}

While \ref{enum3:i} and \ref{enum3:ii} yield explicit bounds, the
combinatorial approach does this for \ref{enum3:iii}.

\begin{theorem}[Explicit approximation]
  \label{thm:E_complete}
  Let $r,q \geq 2$, and $\epsilon_{r,n}$ and $\kappa$ as in
  \autoref{thm:E_by_gen}, and $n$ be composite. Then for the number
  $\# \Eone{r,n}(\Fq)$ of relatively irreducible monic $r$-variate
polynomials of degree $n$ over $\Fq$ we have
    \begin{equation}
    \left| \# \Eone{r,n}(\FF_q) - \epsilon_{r,n}(q) \right| \leq
    \epsilon_{r,n}(q) \cdot 3 q^{-\kappa}.
    \end{equation}
\end{theorem}

\subsection{Reducible space curves}
\label{sec:reducible-curves}

The \emph{Chow variety} of curves of degree $n$ in the $r$-dimensional
projective space $\PP^{r} = \PP^{r}(\Fqbar)$ over an algebraic closure
$\Fqbar$ is
denoted by $C_{r,n}$.  Each point of the Chow variety $C_{r,n}$
actually corresponds to a unique \emph{effective cycle} in $\PP^{r}$ of
dimension $1$ and degree $n$, that is, to a formal linear combination
$\sum a_{i} C_{i}$, where each $C_{i}$ is an irreducible curve in
$\PP^{r}$, each $a_{i}$ is a positive integer and $\sum a_{i} \deg
(C_{i}) = n$.

For a subfield $F\subseteq \Fqbar$, an effective $F$-cycle $C$ is
called \emph{$F$-reducible} if there exist $m \geq 2$ and effective
$F$-cycles $C_{1}, \dots, C_{m}$ such that $C = \sum_{i=1}^{m} C_{i}$
holds. Let $C_{r,n}(\Fq)$ denote the Chow variety of effective
$\Fq$-cycles and $\spaceR{r,n}(\Fq)$ its closed subvariety of
$\Fq$-reducible $\Fq$-cycles. Methods of algebraic geometry yield the
following bounds on the probability that a random curve of degree~$n$ in
$\PP^{r}(\Fq)$ is $\Fq$-reducible.

\begin{theorem}
  \label{thr:1}
Let $r \geq 3$ and
\begin{align}
 g_{r,n} & = {\binom{r+n-2}{n}}^{2} \cdot
 \frac{r+n-1}{(r-1)(n+1)},  \label{eq:34} \\
c_{r,n} & =(2en)^{r(r+1)(n^{2}+1)+4rg_{r,n}},
\end{align}
where $e$ denotes the basis of the natural logarithm. For the
number $\# \spaceR{r,n} (\Fq)$ of $\Fq$-reducible cycles of degree~$n$
we have the
following.
\begin{ronumerate}
    \item If $n \geq \min\{ 4r-7, 7\}$, then
\begin{equation}
  \label{eq:4}
  \frac{1}{4c_{r,n}}q^{-(n-2r+3)} \leq \frac{\# \spaceR{r,n} (\Fq)}{\#
  C_{r,n}(\Fq)} \leq c_{r,n} q^{-(n-2r+3)}.
\end{equation}
\item If $n = 4r-8$, then
\begin{equation}
  \frac{1}{2n!\,c_{r,n}} q^{-r+2} \leq \frac{\# \spaceR{r,n} (\Fq)}{\#
  C_{r,n}(\Fq)} \leq c_{r,n} q^{-r+2}.
\end{equation}
  \end{ronumerate}
\end{theorem}

We call an $\Fqbar$-reducible cycle \emph{absolutely reducible}. An
$\Fq$-cycle can be absolutely reducible for two reasons: either it is
$\Fq$-reducible, as treated above, or \emph{relatively
  $\Fq$-irreducible}, that is, is $\Fq$-irreducible and
$\Fqbar$-reducible.  The set of relatively $\Fq$-irreducible (or
\emph{exceptional}) $\Fq$-curves of degree $n$ in $\PP^{r}$ is denoted by
$\spaceE{r,n} (\Fq)$.

\begin{theorem}
  \label{thr:2}
Let $r \geq 3$, $n \geq 4r - 8$, let $\ell$ denote the smallest
prime divisor of $n$, and
\begin{align}
b_{r,n} & = 3(r-2)+n(n+3)/2, \\
d_{\ell,n,r} & = (en/\ell)^{r(r+1)(n^{2}/\ell^{2}+1)+4rg_{r,n/\ell}}.
\end{align}
For the number $\# \spaceE{r,n}(\Fq)$ of relatively $\Fq$-irreducible
cycles of degree~$n$ we have
\begin{align}
  q^{2n(r-1)}(1-4q^{2(1-n)(r-1)}) \leq \# \spaceE{r,n} (\Fq) \leq
  2d_{\ell,n,r}q^{2n(r-1)} \text{ for } n/\ell \leq 4r-7, \\
q^{\ell b_{r,n/\ell}}(1-16q^{\ell-n}) \leq \# \spaceE{r,n} (\Fq) \leq 3
d_{\ell,n,r}q^{\ell b_{r,n/\ell}} \text{ for } n/\ell \geq 4r-8.
\end{align}
\end{theorem}

\subsection{Decomposable polynomials}
\label{sec:multi-decomp}

For monic univariate $g\in F[y]$ and $h\in \Pone{r,n}$, we define their \emph{composition}
\begin{equation}
  \label{eq:9}
f = g \circ h = g(h) \in \Pone{r,n}.
\end{equation}
If $\deg g \geq 2$ and $\deg h \geq 1$, then $(g,h)$ is a
\emph{decomposition} of $f$. A polynomial $f \in \Pone{r,n}$ is
\emph{decomposable} if there exist such $g$ and $h$. There are other
notions of decompositions. The present one is called
uni-multivariate in \cite{gatgut03}. Another one is studied in
\cite{fauper08} for cryptanalytic purposes. In the context of
univariate polynomials $\deg h \geq 2$ is also
required, see Section~\ref{sec:univ-decomp}.

It is sufficient to concentrate on polynomials with vanishing constant
term, see \autoref{sec:notation}, and we denote by $D_{r,n}(F)$ the
set of all decomposable polynomials $f \in \Pone{r,n}(F)$ with
$f(0,\dots,0) = 0$.

\begin{theorem}
  Let $\Fq$ be a finite field with $q$ elements, $r\geq 2$, and
  $\ell$
  the smallest prime divisor of the composite integer $n \geq 2$.  Let
  \begin{align}\label{eq:pm}
    m & = \begin{cases*}
      n & if $r=2$, $n/\ell$ is prime, and $n/\ell \leq 2\ell - 5$, \\
      \ell & otherwise,
    \end{cases*} \\
    \alpha_{r,n}  &=
    q^{\binom{r+n/m}{r}+m-3} \frac{1-q^{-\binom{r-1+n/m}{r-1}}}{1-q^{-1}}, \\
\beta_{r,n} & = \frac{2q^{-\frac{1}{2} \binom{r-1+n/\ell}{r-1} + 1}}{1-q^{-1}}.
  \end{align}
  Then for the number $\# D_{r,n}(\Fq)$ of decomposable monic
  $r$-variate polynomials with vanishing constant term of degree~$n$
  over $\Fq$ we have
  \begin{equation}
    \label{eq:3}
    \left|{\#D_{r,n}(\Fq)}- \alpha_{r,n}\right| \leq
\alpha_{r,n} \cdot \beta_{r,n}.
  \end{equation}
\end{theorem}

\section{Counting univariate decomposable polynomials}
\label{sec:univ-decomp}

The \emph{composition} of two univariate polynomials $g,h \in F[x]$ over a field
$F$ is denoted as $f= g \circ h= g(h)$, and then $(g,h)$ is a
\emph{decomposition} of $f$, and $f$ is \emph{decomposable} if $g$ and
$h$ have degree at least $2$.  In the 1920s, Ritt, Fatou, and Julia
studied structural properties of these decompositions over
$\mathbb{C}$, using analytic methods. Particularly important are two
theorems by Ritt on the uniqueness, in a suitable sense, of
decompositions, the first one for (many) indecomposable components and
the second one for two components, as above.  \cite{eng41} and
\cite{lev42} proved them over arbitrary fields of characteristic zero
using algebraic methods.

The theory was extended to arbitrary characteristic by
\cite{frimac69}, \cite{dorwha74}, \cite{sch82c, sch00c}, \cite{zan93},
and others. Its use in a cryptographic context was suggested by
\cite{cad85}. In computer algebra, the decomposition method of
\cite{barzip85} requires exponential time.  A fundamental dichotomy is
between the \emph{tame case}, where the characteristic $p$ does not
divide $\deg g$, and the \emph{wild case}, where $p$ divides $\deg g$,
see \cite{gat90d,gat90c}.  (\cite{sch00c}, \S~1.5, uses \emph{tame} in
a different sense.)  A breakthrough result of \cite{kozlan89} was
their polynomial-time algorithm to compute tame decompositions; see
also \cite*{gatkoz87}; \cite*{kozlan96}; \cite{gutsev06}, and the
survey articles of \cite{gat02c} and \cite{gutkoz03} with further
references.  Schur's conjecture, as proven by \cite{tur95}, offers a
natural connection between the tame indecomposable polynomials in this
section and certain absolutely irreducible bivariate polynomials, as
studied in Section~\ref{sec:rel_irr}.  More precisely, a tame
polynomial $f$ is indecomposable if $(f(x)-f(y))/(x-y)$ is absolutely
irreducible.  Aside from natural exceptions, the converse is also
true.

In the wild case, considerably less is known, both
mathematically and computationally. \cite{zip91} suggests that the
block decompositions of \cite{lanmil85} for determining subfields of
algebraic number fields can be applied to decomposing rational
functions even in the wild case.  A version of Zippel's algorithm in
\cite{bla14} computes in polynomial time all decompositions of a polynomial that are
minimal in a certain sense.  \cite{avazan03} study ambiguities in the
decomposition of rational functions over $\mathbb{C}$.  On a different but related topic, \cite{ziemue08} found interesting
characterizations for Ritt's First Theorem, which deals with complete
decompositions, where all components are indecomposable.

We have seen fairly precise estimates for the number of multivariate
decomposable polynomials in Section~\ref{sec:multi-decomp}.  It is
intuitively clear that the univariate decomposable polynomials also
form only a small minority among all univariate polynomials over a
field, and this second part of our survey confirms this intuition.
The task is to approximate the number of decomposables over a finite
field, together with a good relative error bound.  One readily obtains
an upper bound. The challenge then is to find an essentially matching
lower bound.

A set of distinct decompositions of $f$ is called a
\emph{collision}. The number of decomposable polynomials of degree $n$
is thus the number of all pairs $(g,h)$ with $\deg g \cdot \deg h = n$
reduced by the ambiguities introduced by collisions.  An important
tool for estimating the number of collisions is Ritt's Second Theorem.
The first algebraic versions of this in positive characteristic $p$
required $p>\deg(g\circ h)$. \cite{zan93} reduced this to the milder
and more natural requirement $g' \neq 0$ for all $g$ in the collision.
His proof works over an algebraic closed field, and Schinzel's
\citeyear{sch00c} monograph adapts it to finite fields. In
Section~\ref{sec:ritt-2}, we provide a precise quantitative version of
Ritt's Second Theorem, by determining exactly the number of such
collisions in the tame case, assuming that $p \nmid n/\ell$, where $n$
is the degree of the composition and $\ell$ is the smallest prime
divisor of $n$. This is based on a unique normal form for the
polynomials occurring in Ritt's Second Theorem.

\cite{gie88b} was the first to consider this counting problem.  He
showed that the decomposable polynomials form an exponentially small
fraction of all univariate polynomials.  General approximations to the
number of univariate decomposable polynomials are shown in
Section~\ref{sec:cdup}.  They come with satisfactory (rapidly
decreasing) relative error bounds except when $p$ divides $n = \deg f$
exactly twice. \cite{zie14} provides an exact count of tame univariate
polynomials. In Section~\ref{sec:prelim}, we determine exactly the
number of decomposable polynomials in one of the difficult wild cases,
namely when $n = p^{2}$.

\cite{zan08} studies a different but related question, namely
compositions $f= g \circ h$ in $ \mathbb{C} [x]$ with a \emph{sparse}
polynomial $f$, having $t$ terms. The degree is not bounded.  He gives
bounds, depending only on $t$, on the degree of $g$ and the number of
terms in $h$.  Furthermore, he gives a parametrization of all such
$f$, $g$, $h$ in terms of varieties (for the coefficients) and
lattices (for the exponents).  \citet{boddeb09} also deal with
counting.

Unless otherwise attributed, the results of Section~\ref{sec:ritt-2}
are from \cite{gat12a}, those of Section~\ref{sec:cdup} from
\cite{gat08c}, and those of Section~\ref{sec:prelim} from
\cite*{blagat13}.

\subsection{Notation}
\label{sec:notation}

A nonzero
polynomial $f\in F[x]$ over a field $F$ of characteristic $p \geq 0$ is \emph{monic} if
its leading coefficient $\operatorname{lc}(f)$ equals $1$.  We call $f$
\emph{original} if its graph contains the origin, that is, $f(0)=0$.
  For $g, h \in F[x]$,
\begin{equation}
  \label{defComp:f}
  f = g \circ h = g(h) \in F[x]
\end{equation}
is their \emph{composition}.  If $\deg g, \deg h \geq 2$, then $(g,h)$
is a \emph{decomposition} of $f$. A polynomial $f \in F[x]$ of degree
at least $2$ is \emph{decomposable} if there exist such $g$ and $h$,
otherwise $f$ is \emph{indecomposable}.  A decomposition
\eqref{defComp:f} is \emph{tame} if $p\nmid \deg g$, and $f$ is
\emph{tame} if $p\nmid \deg f$.

  Multiplication by a unit or addition of a constant does not change
  decomposability, since
  \begin{equation}
    \label{eq:11}
  f = g \circ h \Longleftrightarrow a f+b = (a g+b) \circ h
  \end{equation}
  for all $f$, $g$, $h$ as above and $a,b \in F$ with $a\neq 0$.  In
  other words, the set of decomposable polynomials is invariant under
  this action of $F^{\times} \times F$ on $F[x]$.

Furthermore, any decomposition $(g,h)$ can be normalized by this
action, by taking $a = \operatorname{lc} (h)^{-1} \in F^{\times}$, $b=-a \cdot h(0)
\in F$, $g^{*} = g((x-b)a^{-1}) \in F[x]$, and $h^{*} = ah+b$.  Then
$g\circ h = g^{*} \circ h^{*}$ and $g^{*}$ and $ h^{*}$ are monic original.

It is therefore sufficient to consider compositions $f = g \circ h$
where all three polynomials are monic original. In such a tame
decomposition, $g$ and $h$ are uniquely determined by $f$ and $\deg
g$. For $n \geq 1$ and any proper divisor $e$ of $n$, we write
\begin{align}
    P_{n}(F) & = \{ f \in F[x] \colon \text{$f$ is monic original of degree $n$}\}, \label{eq:45} \\
    D_{n}(F) & = \{ f \in P_{n}(F) \colon \text{$f$ is
      decomposable} \}, \label{eq:46} \\
    D_{n,e}(F) & = \{ f \in P_{n}(F) \colon \text{$f = g \circ h$
      for some $(g,h) \in P_{e}(F) \times P_{n/e}(F)$} \}.
\end{align}
Thus $P_{n}(F)$ and $D_{n}(F)$ are the subsets of original
polynomials in the sets $P_{1,n}(F)$ and $D_{1,n}(F)$,
respectively, as defined in the context of multivariate polynomials
(\autoref{sec:multi-decomp}) but with right component $h$ of degree at least
$2$.  We sometimes leave out $F$ from the notation when it is clear
from the context and have over a finite field $\mathbb{F}_{q}$ with
$q$ elements
\begin{align}
  \#P_{n} &= q^{n-1},\\
  \#D_{n,e} &\leq q^{e+n/e-2}.
\end{align}
The set $D_{n}$ of all decomposable polynomials in $P_{n}$
satisfies
\begin{equation}\label{substack}
  D_{n}= \bigcup_{\substack{e\mid n\\1<e<n}} D_{n,e}.
\end{equation}
In particular, $D_{n} = \varnothing$ if $n$ is prime and $x \in P_{1}$
is neither decomposable nor indecomposable. For the resulting
inclusion-exclusion formula for $\# D_{n}$, we have to determine the
\emph{collisions} (or nonuniqueness) of decompositions, that is,
different components $(g,h)\neq(g^{*},h^{*})$ with equal composition
$g\circ h=g^{*}\circ h^{*}$.

It is useful to single out a special case of wild compositions when $p
> 0$.
\begin{example}\label{rem:coll}
  We call an $f \in P_{n} \cap F[x^{p}]$ a \emph{Frobenius
    composition}, since then $f = g^{*} \circ x^{p}$ for some $g^{*}
  \in P_{n/p}$, and any decomposition $(g,h)$ of $f = g \circ h$ is a
  \emph{Frobenius decomposition}. We denote by $\varphi \colon F
  \longrightarrow F$ the Frobenius endomorphism over a field $F$ of
  characteristic $p$, with $\varphi(a)= a^{p}$ for all $a \in F$, and
  extend it to an $\mathbb{F}_{p}$-linear map $\varphi \colon P_{n}
  \longrightarrow P_{n}$ with $\varphi(x)=x$. For $h \in P_{n/p}
  \mysetminus \{x^{p}\}$, this provides the collision
  \begin{equation}\label{eq:frob}
    x^{p} \circ h = \varphi (h) \circ x^{p}.
  \end{equation}
\end{example}

If $F$ is perfect -- in particular if $F$ is finite or algebraically
closed -- then $\varphi$ is an automorphism on $F$ and every Frobenius
composition except $x^{p^{2}}$ is a collision as in
\eqref{eq:frob}. Over $F = \mathbb{F}_{q}$, this yields $q^{p - 1} -
1$ collisions in $D_{p^{2}}$ and $q^{n/p-1}$ collisions in $D_{n}$ for
$p \mid n \neq p^{2}$, called \emph{Frobenius collisions}. This example is noted in \citet[Section~I.5, page~39]{sch82c}.

For $f \in P_{n} (F)$ and $a \in F$, the \emph{original shift} of
$f$ by $a$ is
\begin{equation}
  \label{eq:14}
  f^{[a]} = (x-f(a)) \circ f \circ (x+a) \in P_{n}(F).
\end{equation}
Original shifting defines a group action of the additive group of
$F$ on $P_{n}(F)$.  Shifting respects decompositions in the sense
that for each decomposition $(g,h)$ of $f$ we have a decomposition
$(g^{[h(a)]}, h^{[a]})$ of $f^{[a]}$, and vice versa.  We denote
$(g^{[h(a)]}, h^{[a]})$ as $(g, h)^{[a]}$.

\subsection{Normal form for Ritt's Second Theorem}
\label{sec:ritt-2}

Ritt presented two types of essential collisions:
\begin{align}
x^{\ell}\circ x^{k}w(x^{\ell}) & =
x^{k\ell}w^{\ell}(x^{\ell})=x^{k}w^{\ell}\circ
x^{\ell}, \label{al:colli} \\
T_{m}(x,z^{\ell})\circ T_{\ell}(x,z) & = T_{\ell m}(x,z)=T_{\ell}(x,z^{m})\circ
T_{m}(x,z),
\end{align}
where $w\in F[x]$, $z\in F^{\times} = F \mysetminus \{0\}$, and $T_{m}$
is the $m$th Dickson polynomial of the first kind. And then he proved
that these are all possibilities up to composition with linear polynomials.
This involved four unspecified linear functions, and it is not clear whether there is a relation
between the first and the second type of example.

\Citet{gat12a} presents a normal form for the decompositions in Ritt's
Theorem under Zannier's assumption $g'(g^{*})'\neq 0$ and the standard
assumption $\gcd(\ell,m)=1$, where $m=k+\ell\deg w$ in
\eqref{al:colli}. This normal form is unique unless $p\mid m$.

\begin{theorem}[Ritt's Second Theorem, normal form] \label{th:fifi}
  Let $F$ be a field of characteristic $p \geq 0$, let $m>\ell \geq 2$ be
  integers with $\gcd(\ell, m) = 1$ and $n=\ell m$.  Furthermore, we have monic original $f, g, h,
  g^{*}, h^{*} \in F[x]$ satisfying
  \begin{gather}
    f = g \circ h = g^{*} \circ h^{*}, \label{eq:2} \\
    f, g, h, g^{*}, h^{*} \text{ are monic original}, \label{eq:6} \\
    \deg g = \deg h^{*} = m, \deg h = \deg g^{*} = \ell, \label{eq:8} \\
    g'(g^{*})' \neq 0, \label{eq:15}
  \end{gather}
 where $g'= \partial g/ \partial x$ is the derivative of $g$.  Then either \ref{th:fifi-1} or \ref{th:fifi-2} hold, and \ref{th:fifi-3} is also
  valid.   \noeqref{eq:6,eq:8}
  \begin{ronumerate}
  \item\label{th:fifi-1} (First Case) There exists a monic
    polynomial $w \in F[x]$ of degree $s$ and $a \in F$ so that
    \begin{equation}\label{eq:mopo}
      f= (x^{k\ell}w^{\ell}(x^{\ell}))^{[a]} ,
    \end{equation}
   where $m=s\ell+k$ is the division with remainder of $m$ by $\ell$, with
   $1 \leq k < \ell$. Furthermore, we have
    \begin{gather}
    \begin{split}
      (g,h) & = (x^{k}w^{\ell},x^{\ell})^{[a]}, \\
      (g^{*}, h ^{*}) & = (x^{\ell}, x^{k} w (x^{\ell}))^{[a]} ,
    \end{split} \label{eq:unidet-1} \\
    kw+\ell xw' \neq 0 \text{ and } p \nmid \ell. \label{eq:unidet}
    \end{gather}

    Conversely, any $(w,a)$ as above for which \eqref{eq:unidet}
    holds yields a collision satisfying \eqref{eq:2} through
    \eqref{eq:15}, via \eqref{eq:unidet-1}. If $p\nmid m$, then
    $(w,a)$ is uniquely determined by $f$ and $\ell$.
  \item\label{th:fifi-2} (Second Case) There exist $z,a \in F$
    with $z \neq 0$ so that
    \begin{equation}\label{eq:TN}
      f = T_{n}(x,z)^{[a]} .
    \end{equation}
    Now $(z,a)$ is uniquely determined by $f$. Furthermore, we have
    \begin{gather}
      \begin{split}
      (g,h) & = (T_{m}(x,z^{\ell}) ,
      T_{\ell}(x,z))^{[a]}, \\
      (g^{*}, h^{*}) & = (T_{\ell}(x,z^{m}),
      T_{m}(x,z))^{[a]},
    \end{split} \label{eq:ab} \\
     p \nmid n. \label{eq:ab-5}
    \end{gather}
    Conversely, if \eqref{eq:ab-5} holds, then any $(z,a)$ as
    above yields a collision satisfying \eqref{eq:2} through
    \eqref{eq:15}, via \eqref{eq:ab}.

  \item\label{th:fifi-3} When $\ell \geq 3$, the First and Second Cases
    are mutually exclusive. For $\ell=2$, the Second Case is included in
    the First Case.
  \end{ronumerate}
\end{theorem}

If $p \nmid n$, then the case where $\gcd (\ell,m) \neq 1$ is reduced
to the previous one by a result of \cite{tor88a}. This
determines $D_{n,\ell} \cap D_{n,m}$ exactly if $p \nmid n =
\ell m$.
\begin{theorem}[Tame case]
  \label{thm:FFC}
  Let $\Fq$ be a finite field of characteristic $p$, let $\delta$ denote
  Kronecker's delta function, and
  let $m > \ell \geq 2$ be integers with $p \nmid n = \ell m$, $i =
  \gcd(\ell,m)$ and $s=\lfloor m/\ell \rfloor$. For the number of
  monic original polynomials of degree $n$ over $\Fq$ with left
  components of degree $\ell$ and $m$ we have
  \begin{equation}
    \#(D_{n,\ell}(\Fq) \cap
    D_{n,m}(\Fq)) = \begin{cases}
      q^{2\ell+s-3} & \text{if } \ell \mid m,\\
      q^{2i}(q^{s-1}+(1- \delta_{\ell,2})
      (1-q^{-1})) & \\
      \quad \quad \quad \leq q^{2\ell+s-3} & \text{otherwise}.
    \end{cases}
  \end{equation}
\end{theorem}

In the remaining case where $p \mid n$, the Frobenius collisions are
easily counted and therefore excluded. We have the following upper bounds.

\begin{corollary}[Wild case, upper bounds]
  \label{cor:ffchar}
  Let $\Fq$ be a finite field of characteristic $p$ and
  $\ell$, $m$, $n \geq 2$ be integers with $p \mid n = \ell m$, and
  let $c =\#(D_{n, \ell}(\Fq) \cap D_{n,m} (\Fq) \mysetminus F[x^{p}])$ be the
  number of monic original polynomials of degree $n$
  over $\Fq$ with left components of degree $\ell$ and $m$ that are
  not Frobenius collisions. Then the
  following hold.
  \begin{ronumerate}
  \item\label{cor:ffchar-2} If $p \nmid \ell$, then
    \begin{equation}
    c \leq q^{m+\lceil \ell/p \rceil-2}.
    \end{equation}
  \item
    \label{cor:ffchar-3}
    If $p\mid \ell$ and $\ell < m$, we set $b=\lceil (m-\ell+1)/\ell\rceil$. Then
    \begin{equation}
    c\leq q^{m+\ell-b+\lceil b/p\rceil-2}.
    \end{equation}
 \end{ronumerate}
\end{corollary}

For perspective, we also note the following lower bounds on $c$ from
\cite{gat08c,gat13}. Unlike the exact result of \autoref{thm:FFC}, there is a
substantial gap between the upper and lower bounds.

\begin{corollary}[Wild case, lower bounds]
  \label{cor:ffcharb}
  Let $\Fq$ be a finite field of characteristic $p$, $\ell$ a prime
  number dividing $m >\ell$, assume that $p \mid n = \ell m$, and let
  $c =\#(D_{n, \ell} (\Fq) \cap D_{n,m} (\Fq) \mysetminus F[x^{p}])$
  be the number of monic original polynomials of degree $n$ over $\Fq$
  with left components of degree $\ell$ and $m$ that are not Frobenius
  collisions. Then the following hold.
  \begin{ronumerate}
  \item\label{cor:ffcharb-4} If $p=\ell \mid m$ and each nontrivial
    divisor of $m / p$ is larger than $p$, then
    \begin{equation}
    c \geq  q^{2p+m/p-3}(1-q^{-1})(1-q^{-p+1}).
    \end{equation}
  \item
    \label{cor:ffcharb-1}
    If $p \neq \ell$ divides $m$ exactly $d \geq 1$ times, then
    \begin{equation}\label{eq:ffcharb-1}
    c \geq q^{2\ell+m/\ell-3}(1-q^{-m/\ell})(1-q^{-1}(1+q^{-p+2}
    \frac{(1-q^{-1})^{2}}{1-q^{-p}}))
    \end{equation}
    if $\ell \nmid p^{d}-1$. Otherwise we set $\mu=\gcd(p^{d}-1,\ell)$,
    $r=(p^{d}-1)/\mu$ and have
\begin{align}\label{al:ffcharb-1}
  \begin{aligned}
     c & \geq q^{2\ell+m/\ell-3}\bigl((1-q^{-1}(1+q^{-p+2}
      \frac{(1-q^{-1})^{2}}{1-q^{-p}}))(1-q^{-m/\ell})\\
      & \quad-q^{-m/\ell-r+2}
      \frac{(1-q^{-1})^{2}(1-q^{-r(\mu-1)})}{1-q^{-r}}
      (1+q^{-r(p-2)})\bigr).
    \end{aligned}
\end{align}
 \end{ronumerate}
\end{corollary}

\subsection{The number of decomposable univariate polynomials}
\label{sec:cdup}

The basic statement is that $\alpha_{n}$ as in \eqref{eq:17} is an
approximation to the number of monic original decomposable polynomials of degree
$n$, with relative error bounds of varying quality.
The following is a condensed version of the more precise bounds in \cite{gat08c}.
\begin{theorem}
\label{cor:Fq}
  Let $\Fq$ be a finite field with $q$ elements and
  characteristic $p$, let $\ell$ be the smallest prime divisor of the
  composite integer $n \geq 2$, and
  \begin{equation}\label{eq:17}
    \alpha_{n} =
    \begin{cases}
      2q^{\ell+n/\ell-2}
      & \text{if } n \neq \ell^{2}, \\
      q^{2\ell-2}
      & \text{if } n = \ell^{2}.\\
    \end{cases}
  \end{equation}
  Then the following hold for the number $\# D_{n} (\Fq)$ of decomposable
  monic original polynomials of degree $n$ over $\Fq$, where $p \parallel n$ means that $p$ divides $n$ exactly
    twice.
  \begin{ronumerate}
  \item\label{cor:Fq-6} $ q^{2\sqrt{n}-2} \leq \alpha_{n}
    \leq 2q^{n/2}$.
  \item\label{cor:Fq-1} $ \alpha_{n}/2 \leq \# D_{n} (\Fq) \leq
    \alpha_{n}(1+q^{-n/3\ell^{2}}) < 2\alpha_{n} \leq 4q^{n/2}$.
  \item\label{cor:Fq-2} If $n\neq p^{2}$ and $q> 5$, then $\# D_{n} (\Fq) \geq
    (3-2q^{-1})\alpha_{n}/4 \geq q^{2\sqrt{n}-2}/2$.
  \item\label{cor:Fq-4} Unless $p = \ell \parallel n$ and , we have
    $\#D_{n} (\Fq) \geq \alpha_{n}(1-2q^{-1})$.
  \item\label{cor:Fq-5} If $p \nmid n$, then $| \# D_{n} (\Fq) -
    \alpha_{n} | \leq \alpha_{n}\cdot q^{-n/3\ell^{2}}$.
  \end{ronumerate}
\end{theorem}
The relative error in \ref{cor:Fq-5} is exponentially decreasing in
the input size $n \log q$, in the tame case and for growing
$n/3\ell^{2}$. In \ref{cor:Fq-4}, the factor is $1 + O(q^{-1})$ over
$\Fq$. When $p = \ell \parallel n$, then we have a
factor of about $2$ in \ref{cor:Fq-1}, which is improved to about
$4/3$ in \ref{cor:Fq-2}. The case $n = p^{2}$ is settled in
\autoref{sec:prelim}.

Beyond the previous precise bounds, without
asymptotics or unspecified constants, we now derive some conclusions
about the asymptotic behavior.  There are two parameters: the field
size $q$ and the degree $n$. When $n$ is prime, then $\#
D_{n} (\Fq) = 0$, and prime values of $n$ are excepted in the
following. We consider the asymptotics in one parameter, where the
other one is fixed, and also the special situations where
$\gcd(q,n)=1$. Furthermore, we denote as ``$q,n\longrightarrow
\infty $'' any infinite sequence of pairwise distinct
$(q,n)$. The cases $n=4$ and $p^{2}\parallel n \neq p^2$ for some
prime $p$ are the only ones where
our methods do not show that $\# D_{n} (\Fq)/\alpha_{n}\longrightarrow 1$.
\begin{theorem}\label{thm:consid}
  Let $\# D_{n}(\Fq)$ be the number of decomposable monic original
  polynomials of degree $n$ over $\Fq$, $\alpha_{n}$ as in
  \eqref{eq:17}, and $\nu_{q,n} = \# D_{n}(\Fq)/\alpha_{n}$. We only
  consider composite $n$.
  \begin{ronumerate}
  \item\label{thm:consid-1} For any $q$, we have
    \begin{equation}
      \label{eq:25}
    \lim_{\substack{n\to\infty \\ \gcd(q,n)=1}}{\nu_{q,n}}=1,
    \end{equation}
    \begin{equation}
      \label{eq:24}
    \underset{n\to\infty}{\lim\sup} ~ {\nu_{q,n}}=1,
    \end{equation}
    \begin{align}
      \frac 1 2 &\leq \nu_{q,n} \text{ for any }
        n,\\
      \frac{3-2q^{-1}}{4}&\leq \nu_{q,n} \text{ for any }n \text{ if }
        q > 5.
    \end{align}
  \item\label{thm:consid-2} Let $n$ be a composite integer and $\ell$
    its smallest prime divisor. Then
    \begin{equation}
      \label{eq:30}
    \lim_{\substack{q\to\infty\\ \gcd(q,n)=1}}{\nu_{q,n}}=1,
    \end{equation}
    \begin{equation}
      \label{eq:29}
    \underset{q\to\infty}{\lim\sup} ~{\nu_{q,n}}=1,
    \end{equation}
    \begin{equation}
      \underset{{q\to\infty}}{\lim\inf}
      ~{\nu_{q,n}}\begin{cases}
        = 2/3 & \text{ if } n = 4,\\
        \geq\frac 1 4 (3+\frac{1}{\ell+1})\geq \frac 5 6 &\text{ if }
          \ell^{2}\parallel n \text{ and }n\neq \ell^{2},\\
        =1 &\text{ otherwise.}
      \end{cases}
    \end{equation}
  \item\label{thm:consid-3} For any sequence $q,n \rightarrow
    \infty$, we have
    \begin{equation}
      \label{eq:32}
    \lim_{\substack{q,n\to\infty\\ \gcd(q,n)=1}}{\nu_{q,n}}=1,
    \end{equation}
    \begin{equation}
      \label{eq:31}
    \frac 1 2 \leq \underset{q,n\to\infty}{\lim\inf} ~{\nu_{q,n}}\leq
    \underset{q,n\to\infty}{\lim\sup}
    ~{\nu_{q,n}}=1.
    \end{equation}
  \end{ronumerate}
\end{theorem}

\subsection{Collisions at degree $p^{2}$}
\label{sec:prelim}

The previous section gives satisfactory estimates for the number of
decomposable polynomials at degree $n$ unless $p^{2} \parallel n$.
The material of this section determines the number in the easiest
of these open cases, namely for $n = p^{2}$.

First, we present two classes of explicit collisions at degree $r^{2}$, where $r$
is a power of the characteristic $p>0$ of the field $F$.
The collisions of \autoref{thm:nonadd} consist of additive and subadditive polynomials.
A polynomial $\addpol$ of degree $r^{k}$
is \emph{$r$-additive} if it is of the form $\addpol = \sum_{0 \leq i
  \leq k} a_i x^{r^i}$ with all $a_i \in F$. We call a
polynomial \emph{additive} if it is $p$-additive. A polynomial is
additive if and only if it acts additively on an algebraic closure
$\overline{F}$ of $F$, that is $\addpol(a + b) = \addpol(a)+
\addpol(b)$ for all $a$, $b \in \overline{F}$; see \citet[Corollary
1.1.6]{gos96}. The composition of additive polynomials is additive,
see for instance Proposition 1.1.2 of the cited book.  The
decomposition structure of additive polynomials was first studied by
\cite{ore33b}.  \citet[Theorem~4]{dorwha74} show that all components
of an additive polynomial are additive.  \cite{gie88b} gives lower bounds
on the number of decompositions and algorithms to determine them.

For a divisor $m$ of $r -1$, the \emph{$(r,m)$-subadditive} polynomial associated with the $r$-additive polynomial $\addpol$ is $\subadd = x(\sum_{0 \leq i \leq k} a_i x^{(r^i - 1)/m})^m$ of degree $r^{k}$.
Then $\addpol$ and $\subadd$ are related as $x^m \circ \addpol = \subadd \circ x^m$.
\citet{dic97} notes a special case of subadditive polynomials, and \cite{coh85} is concerned with the reducibility of some related polynomials.
\cite{coh90b, coh90c} investigates their connection to exceptional
polynomials and coins the term ``sub-linearized''; see also \cite{cohmat94}.
\citet*{couhav04} derive the number of indecomposable subadditive
 polynomials and present an algorithm to decompose subadditive polynomials.

 \citet[Theorem 3]{ore33b} describes exactly the right components of
 degree $p$ of an additive polynomial.  \cite{henmat99} relate such
 additive decompositions to subadditive polynomials, and in their
 Theorems~3.4 and 3.8 describe the collisions of \autoref{thm:nonadd}
 below.  \autoref{thm:normal} shows that together with
 those of \autoref{thm:constmulti} and the
 Frobenius collisions of \autoref{rem:coll}, these examples and their
 shifts comprise all collisions at degree $p^{2}$.

\begin{fact}
\label{thm:nonadd}
Let $r$ be a power of $p$, $u, s\in F^{\times}$, $\varepsilon \in \{0,1\}$, $m$ a positive divisor of $r-1$, $\ell = (r-1)/m$, and
\begin{equation}
 \label{eq:7}
\begin{split}
  f &= \SimpConst{u,s,\varepsilon,m} = x (x^{\ell(r+1)}-\varepsilon u
  s^{r}x^{\ell} + us^{r+1})^{m} \in P_{r^{2}}(F), \\
 T &= \{t \in F \colon t^{r+1} -\varepsilon ut + u = 0\}.
  \end{split}
\end{equation}
For each $t \in T$ and
\begin{equation}
\label{eq:80}
\begin{split}
g   & = x (x^{\ell}-us^{r}t^{-1})^{m}, \\
h   & = x (x^{\ell}-st)^{m},
\end{split}
\end{equation}
both in $P_{r}(F)$, we have $f = g \circ h$.
Moreover, $f$ has a $\# T$-collision.
\end{fact}

The polynomials $f$ in
\eqref{eq:7} are ``simply original'' in the sense that they
have a simple root at $0$.  This motivates the designation $S$.
The second construction of collisions goes as follows.
\begin{theorem}
\label{thm:constmulti}
  Let $r$ be a power of $p$, $b \in F^\times$, $a \in F \mysetminus \{0, b^{r}\}$, $a^* = b^r - a$, $m$ an integer with $1 < m <r-1$ and $p \nmid m$,  $m^* = r - m$, and
  \begin{equation}
  \label{eq:3normal}
  \begin{split}
   f = \MultConst{a, b, m} &= x^{m m^*} (x - b)^{m m^*} \left(x^m + a^* b^{-r} ((x-b)^{m} - x^m)\right)^m \\
	&\quad\quad\quad \cdot
        \left(x^{m^*} + a b^{-r} ((x-b)^{m^*} - x^{m^*})
        \right)^{m^*},\\
g &= x^m (x - a)^{m^*}, \\
h &= x^{r} + a^{*}b^{-r}(x^{m^*}(x-b)^m - x^{r}),  \\
g^* &= x^{m^*} (x - a^*)^m, \\
h^* &= x^{r} + ab^{-r}(x^m (x-b)^{m^*} - x^{r}).
\end{split}
\end{equation}
Then $f = g \circ h = g^{*} \circ h^{*} \in P_{r^{2}}(F)$ has a 2-collision.
\end{theorem}

The polynomials $f$ in
\eqref{eq:3normal} are ``multiply original'' in the sense that they
have a multiple root at $0$.  This motivates the designation $M$.
The notation is set up so that ${}^{*}$ acts as an involution on our
data, leaving $b$, $f$, $r$, and $x$ invariant.

\cite{zie11} points out that the rational functions of case (4) in Proposition 5.6 of \cite{avazan03} can be transformed into \eqref{eq:3normal}.  Zieve also mentions that this example already occurs in unpublished work of his, joint with Robert Beals.

\begin{theorem}
\label{thm:normal}
Let $F$ be a perfect field of characteristic $p$ and $f \in P_{p^{2}}(F)$. Then $f$ has a collision $\{(g,h), (g^*, h^*)\}$ if and only if exactly one of the following holds.
  \begin{itemize}
    \item[\textnormal{(F)}] \namedlabel{class:0normal}{(F)} The polynomial $f$ is a Frobenius collision as in \autoref{rem:coll}.
    \item[\textnormal{(S)}] \namedlabel{class:1normal}{(S)} The polynomial $f$ is simply original and there are  $u$, $s$,  $\varepsilon$, and  $m$ as in \autoref{thm:nonadd} and $w \in F$  such that
\begin{equation}
f^{[w]} = \SimpConst{u,s,\varepsilon,m}
\end{equation}
and the collision $\{(g,h)^{[w]}, (g^{*},h^{*})^{[w]}\}$ is contained in the collision described in \autoref{thm:nonadd}, with $\# T \geq 2$.
  \item[\textnormal{(M)}] \namedlabel{class:3normal}{(M)} The polynomial $f$ is multiply original and there are  $a$,  $b$, and $m$ as in \autoref{thm:constmulti} and $w \in F$ such that
  \begin{equation}
  f^{[w]} = \MultConst{a,b,m}
  \end{equation}
and the collision $\{(g,h)^{[w]}, (g^*,h^*)^{[w]}\}$ is as in \autoref{thm:constmulti}.
 \end{itemize}
\end{theorem}

In particular, the collisions in case \ref{class:1normal} and case
\ref{class:3normal} have exactly $\# T$ and $2$ distinct decompositions, respectively.  Inclusion-exclusion now yields the following exact formula for the number of decomposable polynomials of degree $p^{2}$ over $\Fq$.

\begin{theorem}
\label{cor:main}
Let $\Fq$ be a finite field of characteristic $p$, $\delta$
Kronecker's delta function, and $\tau$ the number of positive divisors
of $p-1$. Then for the number $\# D_{p^{2}} (\Fq)$ of decomposable monic original
polynomials of degree $p^{2}$ over $\Fq$ we have
  \begin{equation}
    \label{eq:1}
    \begin{split}
    \# D_{p^{2}} (\Fq) & = q^{2p-2} -q^{p-1}  + 1  -  \frac{(\tau q  -q +1)(q-1)(qp-p-2)}{2(p+1)} \\
    & \quad - (1- \delta_{p, 2}) \frac{q(q-1)(q-2)(p-3)}{4}.
    \end{split}
  \end{equation}
\end{theorem}

In particular, we have
\begin{align}
\# D_{4} (\Fq) & = q^{2} \cdot \frac{2+q^{-2}}{3} && \text{for $p =
  2$}, \\
\# D_{9} (\Fq) & = q^{4} \bigl( 1 - \frac{3}{8} ( q^{-1} + q^{-2} - q^{-3} - q^{-4})\bigr) && \text{for $p = 3$},\\
\# D_{p^{2}} (\Fq) & = q^{2p-2} \bigl( 1 - q^{-p+1} + O(q^{-2p+5+1/d}) \bigr) &&\text{for $q = p^{d}$ and $p \geq 5$}.
\end{align}

We have two independent parameters $p$ and $d$, and $q= p^d$. For two eventually positive functions $f, g \colon \mathbb{N}^{2} \rightarrow \mathbb{R}$, here $g \in O(f)$ means that there are constants $b$ and $c$ so that $g(p,d) \leq c \cdot f(p,d)$ for all $p$ and $d$ with $p+d \geq b$.
We have the following asymptotics.
\begin{corollary} Let $p \geq 5$, $d \geq 1$, $q= p^d$, and $k \geq
  1$. Then the number $c_{k}$ of decomposable monic original
  polynomials of degree $p^{2}$ over $\Fq$ with exactly $k$
  decompositions is as follows
\begin{align}
  c_{1} & = q^{2p-2} (1 - 2q^{-p+1} + O(q^{-2p+5+1/d})), \\
  c_{2} & = q^{p-1} (1 + O(q^{-p+4+1/d})), \\
  c_{p+1} & = (\tau - 1) q^{3-3/d} \bigl(1 + O(q^{-\max\{2/d,1-1/d\}})\bigr) \\
& \subseteq O\bigl(q^{3-3/d+1/(d \loglog p)} \bigr), \\
c_{k} & = 0 \text{ if } k \notin \{1,2,p+1\}.
\end{align}
\end{corollary}

\autoref{cor:main} leads to $\lim_{q\to\infty} \nu_{q,\ell^2} = 1$ for
any prime $\ell > 2$ in \autoref{thm:consid}~\ref{thm:consid-2}. For
$n=4$, the sequence has no limit, but oscillates between values close to
$\lim\inf_{q\to\infty} \,\nu_{q,4} = 2/3$ and to $\lim\sup_{q\to\infty}
\,\nu_{q,4} = 1$, and these are the only two accumulation points of
the sequence $\nu_{q,4}$.

\section{Open problems}

Further types of multivariate polynomials that are examined from
a counting perspective include singular bivariate ones \citep{gat08-incl-gat07}
and pairs of coprime polynomials \citep{houmul09}. It remains open to
extend the methods of \autoref{sec:exact} to singular multivariate
ones and achieve exponentially decreasing error bounds for
coprime multivariate polynomials.

For univariate decomposable polynomials, the question of good
asymptotics for $\nu_{q,n}$ when $q$ is fixed and $n \to\infty$ is
still open. More work is needed to understand the case where the
characteristic $p$ is the smallest prime divisor of the degree $n$,
divides $n$ exactly twice, and $n \neq p^{2}$. Ritt's Second Theorem
covers distinct-degree collisions, even in the wild case, see
\cite{zan93}; it would be interesting to see a parametrization even
for $p \mid m$ and obtain a similar classification for general
equal-degree collisions.

Finally, this
survey deals with polynomials only and the study of rational functions
with the same methods remains open.

\section{Acknowledgments}

This work was funded by the B-IT Foundation and the Land Nordrhein-Westfalen.

\nocite{grakal03} 

\bibliographystyle{cc2e}
\bibliography{journals,refs,lncs}

\listoftodos

\end{document}